\magnification=\magstephalf
\input amstex
\documentstyle{amsppt}
\pageheight{49pc}
\pagewidth{33pc}
\NoBlackBoxes
\TagsOnRight
\hcorrection{.25in}
 
\def\az{\alpha}  \def\bz{\beta}
    \def\dz{\delta}
    \def\fz{\varphi}
\def\gz{\gamma}  \def\kz{\kappa}
\def\lz{\lambda} \def\mz{\mu}
     
\def\pz{\pi}

\def\rz{\rho}        \def\sz{\sigma}
\def\tz{\tau}        
\def\vz{\varepsilon}

  \def\ooz{\Omega}

\def\cq{\Cal C}   \def\dq{\Cal D}
\def\eq{\Cal E}

\def\llnbr{\allowlinebreak}

\def\qd{\quad}
\def\qqd{\qquad}

\def\lmm{Lemma}

\def\thm{Theorem}
\def\crl{Corollary}

\def\prf{Proof}
\def\xmp{Example}

\def\d{\text{\rm d}}
\def\op{\text{\rm op}}
\def\jj{J(\d x, \d y)}
\def\jja{J^{(\az)}(\d x, \d y)}

\def\einf{\text{\rm ess\;inf}_{\pi}}
\def\esup{\text{\rm ess\;sup}_{\pi}}
\rightheadtext{Cheeger's inequalities and the existence of spectral gap}
 
\topmatter
\title{Cheeger's inequalities for general symmetric forms 
and existence criteria for spectral gap}\endtitle
\rightheadtext{Cheeger's inequalities and existence criteria for spectral gap}
\author{Mu-Fa Chen and Feng-Yu Wang}\endauthor 
\affil{(Beijing Normal University)
}\endaffil
\thanks{Research supported in part by NSFC (No. 19631060),
Qiu Shi Sci. \& Tech.
Found., DPFIHE, MCSEC and MCMCAS.
Research at MSRI is supported in part by NSF grant DMS-9701755.}\endthanks
\address {Department of Mathematics, Beijing Normal University, 
Beijing 100875,
The People's Republic of China.} \endaddress
\abstract {In this paper, some new forms of the Cheeger's inequalities 
are established for general (maybe unbounded) symmetric forms (\thm\;1.1 and
\thm\;1.2), the resulting estimates improve and extend the
ones obtained by Lawler and Sokal (1988) for bounded jump processes. 
Furthermore, some existence criteria for spectral gap of general 
symmetric forms  or general reversible Markov processes are presented
(\thm\;1.4 and \thm\;3.1), based
on the Cheeger's inequalities and a relationship between the spectral 
gap and the
first Dirichlet and Neumann eigenvalues on local region.
}\endabstract
\endtopmatter
 
\def\d{\text{\rm d}}
\def\op{\text{\rm op}}
\def\jj{J(\d x, \d y)}
\def\jja{J^{(\az)}(\d x, \d y)}

\def\einf{\text{\rm ess\;inf}_{\pi}}
\def\esup{\text{\rm ess\;sup}_{\pi}}
\def\mm#1{\mz_{\langle #1\rangle}^c}
 
\document
 
\head{1. Introduction}\endhead
 
The Cheeger's inequalities$^{[1]}$ are well known and widely used in geometric 
analysis, they provide a practical way to estimate the first eigenvalue
of Laplacian in terms of volumes. These inequalities were then
established for bounded jump processes by Lawler and Sokal$^{[7]}$
(in which, a detail comment on the earlier study and references is  
included). The first aim of the paper is to
establish the inequalities for general (maybe unbounded) symmetric
forms.
 
Let $(E, \eq)$ be a measurable space with reference probability measure
$\pi$. Consider the symmetric form $D$ with domain $\dq(D)$:
$$\align
D(f,g)&=\frac{1}{2}\int J(\d x,\d y)(f(x)-f(y)) (g(x)-g(y))
+\int K(\d x)f(x)g(x), \qqd f, g\in \dq(D)\\
\dq (D)&=\{f\in L^2(\pi):\, D(f, f)<\infty \}.\endalign
$$
where $J$ is a symmetric measure: $J(\d x, \d y)=J(\d y, \d x)$.
Without loss of generality, we assume that $J(\{(x,x)\}: x\in E)=0$.
 
We are interested in the following two quantities:
$$\align
&\lz_0=\inf\{D(f,f):\; \pi (f^2)=1\},\tag 1.1\\
&\lz_1=\inf\{D(f,f):\; \pi (f)=0, \;\pi (f^2)=1\}. \tag 1.2
\endalign$$
We remark that in these definitions, the usual condition ``$f\in \dq(D)$''
is not needed since $D(f,f)=\infty$ for all $f\in L^2(\pi)\setminus \dq(D)$.
We even do not assume in some cases the density of $\dq(D)$ in $L^2(\pi)$.
In what follows, whenever $\lz_1$ is considered, the killing measure
$K(\d x)$ is setting to be zero. In this case, we have $\lz_0=0$ and
$\lz_1$ is known as the spectral gap of the symmetric form $(D, \dq(D))$.
 
Define the Cheeger's constants as follows:
$$\align
&h=\inf_{\pi (A)>0}\frac{J(A\times A^c)+K(A)}{\pi (A)}, \tag 1.3\\
&k=\inf_{\pi (A)\in (0, 1)}\frac{J(A\times A^c)}{\pi (A)\pi (A^c)}, \tag 1.4\\
&k'=\inf_{\pi (A)\in (0, 1/2]}\frac{J(A\times A^c)}{\pi (A)}
   =\inf_{\pi (A)\in (0, 1)}\frac{J(A\times A^c)}{\pi (A)\wedge \pi (A^c)}, 
\tag 1.5\endalign$$
where $a\wedge b=\min\{a, b\}$.
Clearly, $k/2\le k'\le k$ and it is easy to see that $k'$ can be varied over
whole $(k/2, k)$. For instance, take $E=\{0, 1\}$, $K=0$, $J(\{i\}\times
\{j\})=1$ for $i\ne j$ and $\pi(0)=p\le 1/2$, $\pi (1)=1-p$. Then $k'/k=1-p$.
 
Recall that for a given reversible jump process, 
we have a $q$-pair $(q(x), q(x, \d y))$:
$q(x, E)$\llnbr$\le q(x)\le \infty$ for all $x\in E$. Throughout the paper, we
assume that $q(x)<\infty$ for all $x\in E$. The reversibility simply means
that the measure $\pi (\d x) q(x, \d y)$ is symmetric, which gives us
automatically a measure $J$. Then, the killing measure is given by
$K(\d x)=\pi (\d x) d(x)$, where $d(x)=q(x)-q(x, E)$ is called the
non-conservative quantity in the context of jump processes. A jump
process is called bounded if $\sup_{x\in E}q(x)<\infty$. In this case
(or more generally, if $\|J(\cdot, E)+K\|_{\op}<\infty$, 
where $\|\cdot\|_{\op}$ denotes the operator norm from 
$L_+^1(\pi):=\{f\in L^1(\pz):\; f\ge 0\}$ to $\Bbb R_+$, then), 
for the corresponding form, we have $\dq (D)=L^2(\pi)$.
For more details, refer to [2].
 
\proclaim{\thm\;(Lawler \& Sokal)} Take $J(\d x, \d y)=\pi (\d x)
q(x, \d y)$ and suppose that
$\|J(\cdot, E)+K/2\|_{\op}\le M<\infty$. Then, we have
$$h\ge \lz_0\ge \frac{h^2}{2M}. \tag 1.6$$
Next, if additionally $K=0$, then
$$k\ge \lz_1\ge \max\bigg\{\frac{\kz k^2}{8M},\,\frac{k^{'2}}{2M}\bigg\}, 
\tag 1.7$$
where 
$$\kz =\inf_{X, Y}\sup_{c\in\Bbb R}
\frac{(\Bbb E |(X+c)^2-(Y+c)^2|)^2}{1+c^2}\ge 1,$$
the infimum is taken over all i.i.d. random variables $X$ and $Y$ with
$\Bbb E X=0$ and $\Bbb E X^2=1$.\endproclaim
 
In what follows, we consider directly the general symmetric measure $J$ 
whenever it is possible. In other words, we do not require the existence
of a kernel of a modification of $J(\d x, \cdot)/\pi (\d x)$, for which
some extra conditions on $(E, \eq)$ are needed. 
 
We now turn to discuss our general setup. Note that the lower
bounds given in (1.6) and (1.7) decrease to zero as $M\uparrow\infty$. So the  
results would lost their meaning if we go directly from bounded case to 
the unbounded forms. More seriously, when we adopt a general approximation
procedure to reduce the unbounded case to the bounded one 
(cf. [2; \thm\;9.12]), the lower bounds given above usually vanish as we 
go to the limit. To overcome the difficulty, one needs some trick. 
Here we propose a comparison technique. That is, comparing the original form
with some other forms introduced below. 
 
Take and fix a non-negative, symmetric function $r\in\eq\times \eq$ and a 
non-negative function $s\in\eq$ such that
$$\|J^{(1)}(\cdot, E)+ K^{(1)}\|_{\op}\le 1, 
\qqd L_+^1(\pi)\to \Bbb R_+, \tag 1.8$$
where
$$J^{(\az)}(\d x, \d y)=I_{\{r(x, y)>0\}}\frac{\jj}{r(x, y)^{\az}},\qqd
K^{(\az)}(\d x)= I_{\{s(x)>0\}} \frac{K(\d x)}{s(x)^{\az}},
\qqd \az\ge 0. $$
For jump processes, one may simply choose
$$r(x, y)=q(x)\vee q(y)=\max\{q(x), q(y)\}\qd\text{and}\qd s(x)=d(x).$$
We remark that when $\az<1$, the operator 
$J^{(\az)}(\cdot, E) + K^{(\az)}$ from $L_+^1(\pi)$ to $\Bbb R_+$
may no longer be bounded. Correspondingly,
we have symmetric forms $D^{(\az)}$ defined by $(J^{(\az)},\, K^{(\az)})$. 
Therefore, with respect to the form $D^{(\az)}$, 
according to (1.1)---(1.5), we can define $\lz_0^{(\az)}$, $\lz_1^{(\az)}$
and the Cheeger's constants
$h^{(\az)}$, $k^{(\az)}$ and $k^{(\az)'}\,(\az\ge 0)$.
However, in what follows, we need only three cases $\az =0$, $1/2$ and
$1$. When $\az=0$, we return to the original form and so the superscript
``${(\az)}$'' is omitted from our notations.
 
The next two results are our new forms of the Cheeger's inequalities.
  
\proclaim{\thm\;1.1} Suppose that (1.8) holds.
We have
$$\lz_0 \ge \frac{{h^{(1/2)}}^2}{2 -\lz_0^{(1)} }
\ge \frac{{h^{(1/2)}}^2}{1+\sqrt{1-{h^{(1)}}^2}}.\tag 1.9$$
\endproclaim
 
\proclaim{\thm\;1.2} Let $K=0$ and $(1.8)$ hold. Then, we
have
$$\align
\lz_1&\ge \left(\frac{ k^{(1/2)}}{\sqrt{2}+
\sqrt{2-\lz_1^{(1)}}}\,\right)^2,\tag 1.10\\
\lz_1&\ge \frac{{k^{(1/2)}}^{'2}}{1+\sqrt{1-{k^{(1)}}^{'2}}}.\tag 1.11
\endalign$$
\endproclaim
 
When $\|J(\cdot, E)+K\|_{\op}\le M<\infty$, the simplest choice of $r$ 
and $s$ are: $r(x, y)\equiv M$ and $s(x)\equiv M$. Then, (1.8) holds and 
moreover $h^{(1/2)}=h/\sqrt{M}$, ${k^{(1/2)}}'=k'/\sqrt{M}$, 
$h^{(1)}=h/M$ and ${k^{(1)}}'=k'/M$. Hence, by (1.9) and (1.11), we get
$$\lz_0\ge M\big(1-\sqrt{1-h^2/M^2}\,\big)=
\frac{h^2}{M\big(1+\sqrt{1-h^2/M^2}\,\big)}\in 
 \bigg [\frac{h^2}{2M},\,\frac{h^2}{M} \bigg ].$$
and
$$\lz_1\ge M\Big(1-\sqrt{1-k^{'2}/M^2}\,\Big)=
\frac{k^{'2}}{M\big(1+\sqrt{1-k^{'2}/M^2}\,\big)}\in
 \bigg [\frac{k^{'2}}{2M},\, \frac{k^{'2}}{M} \bigg ].\tag 1.12$$
Therefore, for the lower bounds, 
(1.9) improves (1.6) and (1.11) improves (1.7). More essentially, the
lower bound (1.11) is often good enough so that the approximation
procedure [2; \thm\;9.12] mentioned above becomes practical. However,
we will not go to this direction.
In the context of Markov chains on finite graphs, (1.12) was 
obtained before by Chung [5]. Applying (1.12) to $J^{(1)}$, we get
$\lz_1^{(1)}\ge 1-\sqrt{1-{k^{(1)}}^{'2} }$. From this and (1.10), we obtain
$\lz_1\ge\left(\dfrac{ k^{(1/2)} }{\sqrt{2}+
\sqrt{1+\sqrt{1-{k_1^{(1)}}^{'2}} }}\,\right)^2$ which is indeed controlled
by (1.11) since $k^{(\az)}\le 2 k^{(\az)'}$. This means that (1.11) is
usually more practical than (1.10) except a good lower bound of $\lz_1^{(1)}$
is known in advance. However, (1.10) and (1.11) are not comparable even
in the case of $E=\{0, 1\}$. See also the discussion in the second paragraph 
below \lmm\;2.2.
 
In view of \thm\;1.2, we have $\lz_1>0$ whenever $k^{(1/2)}>0$. We now study
some more explicit conditions for the Cheeger's constants appeared in
\thm\;1.2. To state the result, we 
should use the operators corresponding to the forms. 
For a jump process, the operator corresponding to 
$(D^{(\az)}, \dq(D^{(\az)}))$ can be 
expressed by the following simple form
$$\ooz^{(\az)}f(x)=\int I_{[r(x, y)>0]}\frac{q(x, \d y)}{r(x, y)^{\az}}
[f(y)-f(x)]
+ I_{[s(x)>0]}\frac{d(x)}{s(x)^{\az}} f(x).$$ 
 
Next, we need some local quantities of $\lz_0$ and $\lz_1$. 
First, for $B\in\eq$
with $\pi(B)\in (0, 1)$, let 
$\lz_1^{(\az)}(B)$ and $k^{(\az)}(B)$ be defined by (1.2) and (1.4) 
with $E$, $\pi$
and $D$ replaced respectively by $B$, $\pi^B:=\pi(\cdot \cap B)/\pi (B)$ and
$$D_B^{(\az)}(f, f)=\frac{1}{2}\int_{B\times B}\jja (f(y)-f(x))^2.\tag 1.13$$
Second, define
$$\lz_0^{(\az)}(B)=\inf\big\{D^{(\az)}(f,f):\, \pi(f^2)=1,\,f|_{B^c}=0\big\}.$$
As usual, we call $\lz_0^{(\az)}(B)$ and $\lz_1^{(\az)}(B)$ respectively the 
(generalized) first Dirichlet and Neumann eigenvalue on $B$.
It is a simple matter to check that as in (1.7), $k^{(\az)}(B)\ge 
\lz_1^{(\az)} (B)$.
 
For $A\in \eq$, put 
$M_A^{(\az)}=\underset{A}\to{\esup} J^{(\az)}(\d x, A^c)/\pi (\d x)$, where
$\esup$ denotes the essential supremum with respect to $\pi$.
 
\proclaim{\thm\;1.3} Let $K=0$. Given $\az\ge 0$ and $B\in \eq$ with
$\pi (B)>1/2$. Suppose that there exist a function $\fz^{(\az)}$ with 
$\dz_1^{(\az)}\big(\fz^{(\az)}\big):=\text{\rm ess\;sup}_{J^{(\az)}} 
|\fz^{(\az)}(x)-\fz^{(\az)}(y)|<\infty$ and a symmetric operator
$ \big(\ooz^{(\az)}, \dq\big(\ooz^{(\az)}\big)\big)$ 
corresponding to the form 
$\big(D^{(\az)}, \dq\big(D^{(\az)}\big)\big)$ such that
$\dq\big(\ooz^{(\az)}\big)\supset \{I_A:\, A\in\eq,\, A\subset B\}$ and
$\gz^{(\az)}_{B^c}:=-\sup_{B^c}\ooz^{(\az)}\fz^{(\az)}>0 $. Then, 
we have 
$$k^{(\az)}\ge {k^{(\az)}}'\ge
\frac{k^{(\az)}(B)\, \gz_{B^c}^{(\az)}\, [2\pi (B)-1]}
{k^{(\az)}(B) \,\dz_1^{(\az)}\big(\fz^{(\az)}\big)\, [2\pi (B)-1]+ 
\pi (B)^2 \big[\dz_1^{(\az)}\big(\fz^{(\az)}\big) M^{(\az)}_B+ 
\gz_{B^c}^{(\az)}\big]}.$$
\endproclaim
 
Usually, for locally compact $E$, 
we have $k^{(\az)}(B)>0$ and $M_B^{(\az)}<\infty$ for all 
compact $B$. Then the result means that
${k^{(\az)}}'>0$ provided $\dz_1^{(\az)}\big(\fz^{(\az)}\big)<\infty$ and
$\gz_{B^c}^{(\az)}>0$ for large enough $B$. 
 
Up to now, we have discussed the lower bound of $\lz_1$ by using the
Cheeger's constants. However, \thm\;1.3 is indeed a modification of
the second approach we are going to study. That is, 
estimating $\lz_1$ in terms of local $\lz_0$ and $\lz_1$ on subsets
of $E$. The last method has been used recently in the context of diffusions
by Wang [9] and it indeed works for general reversible processes.
The details of the next two results for general situation 
are delayed to Section 3. Here, we 
restrict ourselves to the symmetric forms introduced above.
 
It is the position to state our first criterion for $\lz_1>0$.
 
\proclaim{\thm\;1.4} Let $K=0$. Then for any $A\subset B$ with
$0<\pi(A),\, \pi(B)<1$, we have
$$\frac{\lz_0(A^c)}{\pi(A)}\ge
\lz_1\ge
\frac{\lz_1(B)[\lz_0(A^c)\pi(B)-2M_A\pi(B^c)]}
{2\lz_1(B)+\pi(B)^2 [\lz_0(A^c)+2M_A]}. \tag 1.14$$
\endproclaim
 
As we mentioned before, usually, $\lz_1(B)>0$ for all compact $B$. 
Hence the result means that
$\lz_1>0$ iff $\lz_0(A^c)>0$ for some compact $A$. Because, we can
first fix such an $A$ and then make $B$ large enough so that the right-hand
side of (1.14) becomes positive.
 
Finally, we present an upper bound of $\lz_1$ which provides us a necessary
condition for $\lz_1>0$ and can qualitatively be sharp as illustrated 
by \xmp\;4.5.
 
\proclaim{\thm\;1.5} Let $K=0$, $r>0$, $J$-a.e. and $(1.8)$ hold.
If there exists $\fz\ge 0$ such that 
$$0<\dz_2(\fz):=\text{\rm ess\;sup}_J^{}
|\fz(x)-\fz(y)|^2 r(x, y)<\infty,$$
then
$$\lz_1\le \frac{\dz_2(\fz)}{4}
\sup\Big\{\vz^2:\, \vz\ge 0,\, \pi \big(e^{\vz \fz}\big)
<\infty\Big\}.$$
Consequently, $\lz_1=0$ if there exists $\fz\ge 0$ with 
$0<\delta _2(\fz)<\infty$ such that
$\pi \big(e^{\vz \fz}\big)=\infty$ for all $\vz>0$.
In particular, when $\jj=\pz (\d x)q(x, \d y)$, $\dz_2(\fz)$ can be replaced
by $\dz_2'(\fz):=\esup\int|\fz(x)-\fz(y)|^2 q(x, \d y)<\infty$, without 
using the function $r$ and $(1.8)$.
\endproclaim
 
We mention that the study on the leading eigenvalue of a bounded integral 
operator is indeed included in our general setup. Consider the operator
$P$ on $L^2(\pi)$: $Pf(x)=\int p(x, \d y) f(y)$, generated by an
arbitrary kernel $p(x, \d y)$ with $M:=\sup_x p(x, E)<\infty$. Let 
$\pi (\d x) p(x, \d y)$ be symmetric for a moment. Clearly, the spectrum
of $P$ on $L^2(\pi)$ is determined by the one of $M-P$. Note that
$$\big\langle f,\, (M-P)f\big\rangle_\pz
=\frac{1}{2}\int\pi(\d x) p(x, \d y)\big[f(x)-f(y)\big]^2
+\int\pi(\d x)\big[M- p(x, E)\big] f(x)^2.$$
Thus, the largest (non-trivial) eigenvalue of the integral operator $P$ 
can be deduced from $\lz_0$ or $\lz_1$ treated in the paper. Finally,
by using a symmetrizing procedure, all the results presented here can be 
extended to the non-symmetric forms. Refer to [2; Chapter 9] or [7] 
for instance.
 
The remainder of the paper is organized as follows. Section 2 is devoted
to the proofs of \thm s 1.1---1.3. At the end of the section, a different
approach to handle the unbounded symmetric forms is presented.
A general  existence criterion
for spectral gap is presented in Section 3, which also contains the proofs
of \thm s 1.4 and 1.5. All the results concerning with the spectral gap
are illustrated by Markov chains in the last section.
 
\head{2. Proofs of Theorems 1.1---1.3}\endhead
 
We begin this section with the functional representation of the Cheeger's
constants. The proof is essential the same as in [7] and [8; \S 3.3] for the
bounded situation and hence omitted.	
 
\proclaim{\lmm\;2.1} For every $\az\ge 0$, we have
$$\align
&h^{(\az)}=\inf\bigg\{\frac{1}{2}\int \jja |f(x)-f(y)|+K^{(\az)}(f):\, 
f\ge 0,\, \pi (f)=1 \bigg\}, \\
&k^{(\az)}=\inf\bigg\{\!\int\! \jja |f(x)\!-\!f(y)|:\, f\in L_+^1(\pi),\,
 \int\!\!\pi (\d x)\pi (\d y) |f(x)\!-\!f(y)|\!=\!1\!\bigg\}\\
&\qqd =\inf\bigg\{\int \jja |f(x)-f(y)|:\, f\in L_+^1(\pi),\,
 \pi (|f-\pi (f)|)=1 \bigg\},\\
&k^{(\az)'}=\inf\bigg\{\frac{1}{2}\int \jja |f(x)-f(y)|:\, f\in L_+^1(\pi),\,
\min_{c\in \Bbb R}\pi (|f-c|)=1 \bigg\}. \endalign
$$  \endproclaim
 
\demo{Proof of \thm\;$1.1$} The idea of the proof is based on [7].
 
Let $E^*=E\cup \{\infty\}$. For any $f\in \eq$, define $f^*$ on $E^*$
by setting $f^*=f I_E $. Next, define $J^{*(\az)}$ on $E^*\times E^*$ by
$$J^{*(\az)}(C)=
\cases
J^{(\az)}(C), \qd &C\in \eq\times \eq,\\
K^{(\az)}(A), \qd &C=A\times \{\infty\} \; \text{or} \; \{\infty\}\times A, 
        \; A\in \eq,\\
0, \qd & C=\{\infty\}\times \{\infty\}.
\endcases$$
We have $J^{*(\az)}(\d x, \d y)=J^{*(\az)}(\d y, \d x)$ and
$$\align
\text{\hskip-3em}&\int J^{(\az)}(\d x, E) f(x)^2+ K^{(\az)}(f^2)
=\int J^{*(\az)}(\d x, \d y)f^*(x)^2,
\text{\hskip-2em}\tag 2.1\\
\text{\hskip-3em}&D^{(\az)}(f, f)
=\frac{1}{2}\int J^{*(\az)}(\d x, \d y)(f^*(y)-f^*(x))^2, 
\text{\hskip-2em}\tag 2.2\\
\text{\hskip-3em}&
\frac{1}{2}\int\!\jja |f(y)\!-\!f(x)|\!+\!\!\int\! K^{(\az)}(\d x)|f(x)|
\!=\!\frac{1}{2}\!\int\! J^{*(\az)}(\d x, \d y)|f^*(y)\!-\!f^*(x)|. 
\text{\hskip-2em}\tag 2.3
\endalign$$
Therefore, for $f$ with $\pi (f^2)=1$, by (2.1)--(2.3), (1.8), part (1) of
\lmm\;2.1 and Cauchy-Schwarz inequality,
$$\align
{h^{(1)}}^2&\le
\bigg\{\frac{1}{2}
\int J^{*(1)}(\d x, \d y)\big|f^*(y)^2-f^*(x)^2\big|\bigg\}^2\\
&\le \frac{1}{2} 
D^{(1)} (f, f)\int J^{*(1)}(\d x, \d y)\big[f^*(y)+f^*(x)\big]^2\\
&=\frac{1}{2}
D^{(1)} (f, f)\bigg\{2\int J^{*(1)}(\d x, \d y)
  \big[f^*(y)^2\!+\!f^*(x)^2\big]\!-\!\!
\int J^{*(1)}(\d x, \d y)\big[f^*(y)\!-\!f^*(x)\big]^2\bigg\}\\
&\le D^{(1)} (f, f)\big[2- D^{(1)} (f, f)\big].
\endalign$$
This implies that $D^{(1)}(f, f)\ge 1-\sqrt{1-{h^{(1)}}^2}$ and so
$$\lz_0^{(1)}\ge 1-\sqrt{1-{h^{(1)}}^2 }.\tag 2.4$$
Next, by (1.8), part (1) of \lmm\;2.1 and 
another use of the Cauchy-Schwarz inequality, we obtain
$$\align
{h^{(1/2)}}^2&\le
\bigg\{\frac{1}{2}
\int J^{*(1/2)}(\d x, \d y)\big|f^*(y)^2-f^*(x)^2\big|\bigg\}^2\\
&\le \frac{1}{2}
D (f, f)\int J^{*(1)}(\d x, \d y)\big[f^*(y)+f^*(x)\big]^2\\
&\le D (f, f)\big[2- D^{(1)}(f, f)\big]
\le D (f, f)\big[2- \lz_0^{(1)}\big].  \tag 2.5
\endalign$$
>From this and (2.4), the required assertion follows.\qed
\enddemo
 
\demo{Proof of \thm\;$1.2$} a) First, we prove (1.10). 
Let $f\in\dq(D)$ with $\pi(f)=0$
and $\pi(f^2)=1$. Set $g=f+c$, $c\in\Bbb R$.
Similar to (2.5), we have
$$\align \bigg\{
\int J^{(1/2)}(\d x, \d y)\big|g(y)^2-g(x)^2\big|\bigg\}^2
&\le 4 D (f, f)\big[2(1+c^2)- D^{(1)}(f, f)\big]\\
&\le 4 D (f, f)[2(1+c^2)- \bz]\endalign$$
for all $\bz:\, 0\le \bz<\lz_1^{(1)}\le 2$. 
Hence by \lmm\;2.1, we have
$$D(f, f)\ge \frac{1}{4[ 2(1+c^2)-\bz ]}\bigg\{
\int J^{(1/2)}(\d x, \d y)\big|g(y)^2-g(x)^2\big|\bigg\}^2
\ge \frac{\kz_{\bz}^{}}{4} {k^{(1/2)}}^2 \tag 2.6$$
where $\kz_{\bz}^{}$ is the same as $\kz$ defined below (1.7) but 
replacing the
denominator $1+c^2$ with $2(1+c^2)-\bz$. To estimate $\kz_{\bz}^{}$, we adopt
an optimizing procedure which will be used several times subsequently.  
Set $\gz=\Bbb E |X|\in (0, 1]$. It is known that 
$$\lim_{c\to\pm \infty}\frac{\big(\Bbb E\big|(X+c)^2-(Y+c)^2\big|\big)^2}
{2(1+c^2)-\bz}=2(\Bbb E|X-Y|)^2\ge 2 (\Bbb E |X|)^2=2\gz^2$$
and when $c=0$, $\Bbb E \big|X^2-Y^2\big|\ge 2(1-\Bbb E |X|)=2(1-\gz)$
(cf. [7] or [2; \S 9.2]). Thus, 
$$\kz_{\bz}^{}\ge \inf_{\gz\in (0, 1]}\max\bigg\{2\gz^2,\; \frac{4(1-\gz)^2}
{2-\bz}\bigg\}. \tag 2.7$$
We now need an elementary fact.
 
\proclaim{\lmm\;2.2} Let $f$ and $g$ be continuous functions on $[0, 1]$ and
satisfy $f(0)< g(0)$ and $f(1)>g(1)$. Suppose that $f$ is increasing and
$g$ is decreasing. Then
$$\inf_{\gz\in [0, 1]}\max\{f(\gz),\, g(\gz)\}= f(\gz_0),$$
where $\gz_0$ is the unique solution to the equation $f=g$ on $[0, 1]$.
\endproclaim
 
Applying \lmm\;2.2 to (2.7), we get
$$\kz_{\bz}^{}\ge \frac{4}{\big(\sqrt{2}+\sqrt{2-\bz}\,\big)^2}.$$
Combining this with (2.6) and then letting $\bz\uparrow \lz_1^{(1)}$, we 
obtain (1.10).
 
It is worthy to mention that the estimate just proved can be sharp. To see
this, simply consider $E=\{0, 1\}$, $J(\{i\},\{j\})=1\,(i\ne j)$ and 
$\pi_0=\pi_1=1/2$. Then $k^{(1/2)}=\lz_1^{(1)}=\lz_1=2$.
Moreover, the same example shows that in contract to (1.9), the analog 
of (1.9) 
``$\lz_1\ge  {k^{(1/2)}}^2\big/\big[4\big(2-\lz_1^{(1)}\big)\big]$'' or
``$\lz_1\ge  {k^{(1/2)}}^{'2}\big/\big[2-\lz_1^{(1)}\big)\big]$'' 
does not hold. 
 
b) Define
$$\widetilde D_B^{(\az)}(f, g)=
\frac{1}{2}\int_{B\times B}\jja [f(y)-f(x)]^2
+\int_B J^{(\az)}(\d x, B^c) f(x)^2.$$
It is easy to see that
$$\lz_0(B)=\inf\big\{\widetilde D_B(f, f):\, \pi (f^2 I_B)=1\big\}.$$
Let
$$h_B^{(\az)}=\!\!\!\inf_{A\subset B,\, \pi (A)>0}\!\!
\!\frac{J^{(\az)}(A\!\times\! (B\setminus A))\!+\! 
J^{(\az)}(A\!\times\! B^c)}{\pi (A)}\!
=\!\!\!\inf_{A\subset B,\, \pi (A)>0}\!\!\!
\!\frac{J^{(\az)}(A\!\times\! A^c)}{\pi (A)}.\tag 2.8$$
Then by \thm\;1.1, we have
$\lz_0(B)\ge {h_B^{(1/2)}}^2\big/\Big[1+\sqrt{1-{h_B^{(1)}}^2}\,\Big]$.
 
Next, we prove that
$$\lz_1\ge \inf_{\pi(B)\le 1/2} \lz_0(B).$$
For each $\vz>0$, choose $f_{\vz}$ with $\pi(f_{\vz})=0$ and 
$\pi(f_{\vz}^2)=1$ such that $\lz_1+\vz \ge D(f_{\vz}, f_{\vz})$. Next,
choose $c_{\vz}^{}$ such that $\pi (f_{\vz}<c_{\vz}^{})$, $\pi(f_{\vz}>
c_{\vz}^{})\le 1/2$.
Set $f_{\vz}^{\pm}=(f_{\vz}-c_{\vz}^{})^{\pm}$ and $B_{\vz}^\pm
=\{f_{\vz}^\pm>0\}$.
Then
$$\align
\lz_1+\vz&\ge D(f_{\vz}-c_{\vz}^{},\, f_{\vz}- c_{\vz}^{})\\
&= \frac{1}{2}\int \jj  \big [\big |f_{\vz}^+(y)-f_{\vz}^+(x) \big |+
 \big |f_{\vz}^-(y)-f_{\vz}^-(x) \big |\big]^2\\ 
&\ge  \frac{1}{2}\int \jj  \big(f_{\vz}^+(y)-f_{\vz}^+(x)\big)^2+
  \frac{1}{2}\int \jj  \big(f_{\vz}^-(y)-f_{\vz}^-(x)\big)^2\\
&\ge \lz_0 \big ( B_{\vz}^+\big )\pi \big (  \big (f_{\vz}^+ \big )^2\big )+
   \lz_0 \big ( B_{\vz}^-\big )\pi \big (  \big (f_{\vz}^- \big )^2\big )\\
&\ge \inf_{\pi(B)\le 1/2} 
\lz_0 (B) \pi \big (  \big (f_{\vz}^+ \big )^2+\big (f_{\vz}^- \big )^2\big )
=(1+c_{\vz}^2)\inf_{\pi(B)\le 1/2} \lz_0 (B)
\ge \inf_{\pi(B)\le 1/2} \lz_0 (B).    
\endalign$$
Because $\vz$ is arbitrary, we obtain the required conclusion.
 
Finally, combining the above two assertions together, 
we obtain
$$\align     \lz_1 
&\ge \inf_{\pi(B)\le 1/2} \dfrac{{h_B^{(1/2)}}^2}
{1+\sqrt{1-{h_B^{(1)}}^2}} 
\ge \inf_{\pi(B)\le 1/2} \dfrac{\inf_{\pi(B)\le 1/2}{h_B^{(1/2)}}^2}
{1+\sqrt{1-{h_B^{(1)}}^2}} \\
&\ge \dfrac{\inf_{\pi(B)\le 1/2}{h_B^{(1/2)}}^2}
{1+\sqrt{1-\inf_{\pi(B)\le 1/2}{h_B^{(1)}}^2}}
= \dfrac{{k^{(1/2)}}^{'2}}{1+\sqrt{1-{k^{(1)}}^{'2} }}\qed
\endalign$$\enddemo
 
\demo{\prf\; of \thm\;$1.3$} The proof is split into two lemmas given
below. Noticing that $\az$ is fixed, we may and will omit the superscript 
``$(\az)$'' everywhere in the next two lemmas and their proofs for simplicity. 
\qed\enddemo
 
\proclaim{\lmm\;2.3} Let $B\in\eq$ with $2\pi(B)>1$. Then
$$k'\ge \frac{h_{B^c} k(B)(2\pi (B)-1)}
{k (B)(2\pi (B)-1)+ 2\pi(B)^2(M_B+h_{B^c})},$$
where $h_B$ is defined by (2.8).
\endproclaim
 
\demo{\prf} We need only to consider the case that $h_{B^c} k(B)>0$.
For any $A\in\eq$ with $\pi(A)\in (0, 1/2]$, let $\gz=\pi (AB)/\pi (A)$.
Then
$$\align
\frac{J(A\times A^c)}{\pi(A)}
&= \frac{1}{2\pi (A)}\!\int \jj \big[I_A(y)\!-\!I_A(x)\big]^2
\ge \frac{1}{2\pi (A)}\int_{B\times B}\!\!\jj \big[I_A(y)\!-\!I_A(x)\big]^2\\
&\ge \frac{k(B)\pi^B (A) \pi^B(A^c)}{\pi(A)}
\ge \frac{\pi(B)-1/2}{\pi(B)^2}k(B) \gz. \tag 2.9
\endalign$$
Here, in the last step, we have used $\pi(AB)\le\pi(A)\le 1/2$. 
On the other hand, we have
$$h_{B^c}\pi(AB^c)\le
\frac{1}{2}\int\jj\big[I_{AB^c}(x)-I_{AB^c}(y)\big]^2
=\frac{1}{2}\int\jj\big|I_{A^c\cup B}(x)-I_{A^c\cup B}(y)\big|.
$$
Noticing that $J$ is symmetric and
$$\big|I_{A^c\cup B}(x)-I_{A^c\cup B}(y)\big|
\le \big|I_{A^c}(x)-I_{A^c}(y)\big|+
I_{B\times B^c+B^c\times B}\big|I_{AB}(x)-I_{AB}(y)\big|,$$
we obtain
$$h_{B^c}(1-\gz)=\frac{h_{B^c}\pi(AB^c)}{\pi(A)}
\le \frac{J(A\times A^c)}{\pi(A)}+M_B \gz.$$
Combining this with (2.9) and applying \lmm\;2.2, we get
$$\align
\frac{J(A\times A^c)}{\pi(A)}&\ge
\inf_{\gz\in [0, 1]}\max\Big\{(\pi(B)-1/2)\pi(B)^{-2}k(B) \gz,\;
h_{B^c}-\big(M_B+h_{B^c}\big)\gz \Big\}\\
&=\frac{h_{B^c}k(B)(2\pi (B)-1)}
{k(B)(2\pi (B)-1)+ 2\pi(B)^2(M_B+h_{B^c})}.\qed\endalign
$$\enddemo
 
\proclaim{\lmm\;2.4} Let $\fz$ satisfy $\dz_1(\fz)<\infty$. If 
$\gz_B=-\sup_B\ooz \fz>0$, then $h_B\ge \gz_B/\dz_1(\fz)>0$.
\endproclaim
 
\demo{\prf} For any $A\subset B$, we have
$$\align
\gz_B \pi(A)
&\le \int_A[-\ooz \fz] \d \pi
=\frac{1}{2}\int\jj (I_A(x)-I_A(y))(\fz(x)-\fz(y))\\
&\le \frac{\dz_1(\fz)}{2}\int\jj |I_A(x)-I_A(y)|
=\dz_1(\fz) J(A\times A^c).
\endalign$$
Hence, $h_B\ge \gz_B/\dz_1(\fz)$. \qed\enddemo
 
To conclude this section, we discuss a different way to deal with 
the general symmetric forms. In contrast to the previous approach, we now
keep $(J, K)$ to be the same but change the $L^2$-space. 
To do so, let $p$ be a
measurable function and satisfy $\az_p:=\einf p>0$, $\bz_p:=\pi(p)<\infty$
and $\|J (\cdot, E)+ K \|_{\op}\le \bz_p $\,($L_+^1(\pi_p)\to \Bbb R_+$), 
where $\pi_p=p\pi/\bz_p$. For jump processes, one
may take $p(x)=q(x)\vee r$ for some $r\ge 0$. From this, one sees the
main restriction of the present approach: $\int\pi(\d x) q(x)<\infty$,
since we require that $\pi(p)<\infty$. Except this point, the approach
is not comparable with the previous one (see \xmp\;4.5 and \xmp\;4.7
given below).
 
Next, define $h_p$, $k_p$ and $k_p'$
by (1.3)--(1.5) respectively with $\pi$ replaced by $\pi_p$ and then divided
by $\bz_p$. For instance, 
$k'_p=\inf_{\pi_p(A)\le 1/2} J(A\times A^c)/\pi(pI_A)$.
 
\proclaim{\thm\;2.5} Let $p,\, \az_p,\, \bz_p$ and $\pz_p$ be given above.
Define $\lz_{p, i}\,(i=0, 1)$ by (1.1) and (1.2) with $\pi$ replaced by
$\pi_p$. Then, we have
$$\lz_i \ge \frac{\az_p}{\bz_p}\, \lz_{p, i}, \qqd i=0,\;1.\tag 2.10$$
In particular,
$$\lz_0\ge \az_p\Big(1-\sqrt{1-h_p^2}\,\Big) \tag 2.11$$
and when $K=0$, 
$$\lz_1\ge\max\Bigg\{\frac{\kz}{8} \az_p k_p^2,\;
\az_p\Big(1-\sqrt{1-k_p^{'2}}\,\Big)\bigg\}. \tag 2.12$$
\endproclaim
 
\demo{\prf} a) We prove that $L^\infty(\pz)$ is dense in
$\dq(D)$ in the $D$-norm: $\|f\|_D^2=D(f, f)+\pi(f^2)$.
The proof is similar to [2; \lmm\;9.7]. First, we show that 
$L^\infty(\pz)\subset \dq (d)$. Because $1\in L^1(\pz_p)$ and
$\|J (\cdot, E)+ K \|_{\op}\le \bz_p $, we have 
$J (E, E)+ K(E)\le \bz_p<\infty $. Thus, 
$$\align
D(f, f)&\le \int \jj \big[f(y)^2+f(x)^2\big]+\int K(\d x)f(x)^2\\
&\le 2 \|f\|_\infty^2\big(J(E, E)+K(E)\big)<\infty,
\endalign$$
and hence $f\in\dq (D)$. Next, let $f\in \dq(D)$ and set
$f_n=(-n)\vee (f\wedge n)$. Then $f_n\in \dq(D)$,
$$|f_n(y)-f_n(x)|\le |f(y)-f(x)|\qqd\text{and}\qqd
|f_n(x)|\le |f(x)|\tag 2.13$$
for all $x$, $y$ and $n$. Clearly, $\pz \big((f_n-f)^2\big)\to 0$.
Moreover, since $D(f_n-f, f_n-f)\le 4 D(f, f)<\infty$ by (2.13), we have
$D(f_n-f, f_n-f)\to 0$ by (2.13) and the dominated convergence theorem.
Therefore, $\|f_n-f\|_D\to 0$.
 
b) Here, we prove (2.10) for $i=1$ only since the proof
for $i=0$ is similar and even simpler. Then, (2.11) and (2.12)
follows from (1.7) and the comment right after \thm\;1.2 with 
$M=\bz_p$.
 
Because $L^{\infty}(\pi)\subset L^2(\pz_p)$ and $L^2(\pi_p)$ is just the
domain of the form $D(f, f)$ on $L^2(\pi_p)$, by definition of $\lz_1$
and $\lz_{1, p}$, it suffices to show that
$\pi_p(f^2)-\pi_p (f)^2\ge \big[\pi(f^2)-\pi(f)^2\big]\az_p/\bz_p$
for every $f\in L^\infty(\pi)$. The proof goes as follows.
$$\align
\pi_p\big(f^2\big)-\pi_p(f)^2
&=\inf_{c\in \Bbb R} \int \big(f(x)-c\big)^2 \pi_p(\d x)\\
&=\bz_p^{-1} \inf_{c\in \Bbb R} \int \big(f(x)-c\big)^2 p(x)\pi(\d x)\\
&\ge \frac{\az_p}{\bz_p} \inf_{c\in \Bbb R} \int \big(f(x)-c\big)^2 \pi(\d x)\\
&= \frac{\az_p}{\bz_p}
\big[\pi\big(f^2\big)-\pi(f)^2\big].\qed
\endalign$$\enddemo
 
\head{3. A Criterion for the Existence of Spectral Gap. 
Proofs of Theorem 1.4 and Theorem 1.5}\endhead
 
To state our main criterion, we need some preparation.
 
Let $E$ be a locally compact separable metric space with Borel field
$\eq$ and $\text{supp}(\pi)=E$. Denote by $C_b(E)$ (resp. $C_0(E)$) the set
of all bounded continuous functions (resp. with compact support) on $E$.
 
Next, let $(D, \dq(D))$ be a regular conservative Dirichlet form on
$L^2(\pi)$. By Beurling-Deny's formulae, the form can be expressed as
follows
$$D(f, f)=D^{(c)}(f, f)+
\frac{1}{2}\int\jj (f(x)-f(y))^2,\qqd f\in\dq(D)\cap C_0(E)\tag 3.1$$
where $\dq(D^{(c)})=\dq(D)\cap C_0(E)$ and satisfies a strong local
property; $J$ is a symmetric Radon measure on the product space $E\times E$
off diagonal. Moreover, there exists a finite, non-negative Radon measure 
$\mm f$ such that
$$D^{(c)}(f, f)=\frac{1}{2}\int_E \d \mm f,\qqd f\in\dq(D)\cap C_b(E).$$
 
\proclaim{\thm\;3.1} Let $\cq\subset \dq(D)\cap C_0(E)$ be dense in
$\dq(D)$ in the $D$-norm:
$\|f\|_D^2=D(f, f)+\pi(f^2)$. Set $\cq_L=\{f+c:\, f\in \cq,\; c\in\Bbb R\}$.
Given $A,\, B\in\eq$, $A\subset B$ with $0<\pi(A),\, \pi(B)<1$. Suppose
that the following conditions hold.
\roster
\item There exists a conservative Dirichlet form $(D_B, \dq(D_B))$ on
the square-integrable functions on $B$ with respect to $\pi^B$ such that
$\cq|_B\subset \dq(D_B)$ and
$$D(f, f)\ge D_B(fI_B, fI_B), \qqd f\in \cq_L.$$
\item There exists a function $h\in\cq_L$: $0\le h\le 1$, $h|_A=0$ and 
$h|_{B^c}=1$ such that
$$c(h):=\sup_{f\in \cq_L} \frac{1}{\pi\big(f^2I_B\big)}
\bigg[\frac{1}{2}\int f^2\d \mm h +
\int_{B\times A^c}\jj \big[f(1-h)(y)-f(1-h)(x)\big]^2\bigg]<\infty.$$
\endroster
Then, we have
$$\frac{\lz_0(A^c)}{\pi(A)}\ge
\lz_1\ge
\frac{\lz_1(B)[\lz_0(A^c)\pi(B)-2c(h)\pi(B^c)]}
{2\lz_1(B)+\pi(B)^2[\lz_0(A^c)+2c(h)]}. $$  
\endproclaim
 
\demo{\prf} The upper bound is easy. Simply take $f\in \cq_L$ with $f|_A=0$
and $\pi(f^2)=1$. Then
$$\pi(f^2)-\pi(f)^2=1-\pi\big(f I_{A^c}\big)^2\ge 1-\pi \big(f^2\big)
\pi \big(A^c\big)=1-\pi \big(A^c\big)=\pi (A).$$
Hence $\lz_1\le D(f, f)/\pi(A)$ which gives us 
$\lz_1\le \lz_0\big(A^c\big)/\pi (A)$.
 
For the lower bound, let $f\in \cq_L$ with $\pi(f)=0$ and 
$\pi \big(f^2\big)=1$. Set $\gz =\pi \big(f^2I_B\big)$.
 
a) By condition (1), we have
$$\align
D(f, f)&\ge D_B(fI_B, fI_B)
\ge \lz_1(B)\pi(B)^{-1}\big[\pi \big(f^2I_B\big)-
\pi(B)^{-1}\pi \big(fI_B\big)^2\big]\\
&= \lz_1(B)\pi(B)^{-1}\big[\pi \big(f^2I_B\big)-
\pi(B)^{-1}\pi \big(fI_{B^c}\big)^2\big]\\
&\ge \lz_1(B)\pi(B)^{-1}\big[\gz -
\pi(B)^{-1}\pi \big(f^2I_{B^c}\big)\pi \big(B^c\big)\big]\\
&= \lz_1(B)\pi(B)^{-2}\big[\gz -\pi \big(B^c\big)\big]. \tag 3.2
\endalign$$ 
 
b) Let $\rz$ be the metric in $E$. By the construction of $\mm f$ 
([6; \S 3.2]), there exist a sequence of relatively compact open sets
$G_{\ell}$ increasing to $E$, a sequence of symmetric, non-negative
Radon measures
$\sz_{\bz_n}^{}$ and a sequence $\dz_{\ell}$ such that
$$\align
\int_Eg\d \mm f&=\lim_{\ell\to\infty}\lim_{\bz_n\to\infty}\bz_n
\int_{G_{\ell}\times G_{\ell},\,\rz(x, y)<\dz_{\ell}}
[f(x)-f(y)]^2 g(x)\sz_{\bz_n}^{}(\d x, \d y)\\
&\qqd\qqd\qqd\qqd\qqd\qqd\qqd\qqd f, g\in \dq(D)\cap C_0(E).
\endalign$$
>From this and
$$[(fh)(x)-(fh)(y)]^2 \le 2 h(y)^2 [f(x)-f(y)]^2+
2 f(x)^2 [h(x)-h(y)]^2,$$
it follows that
$$\int \d \mm{fh}
\le 2\int h^2\d \mm f + 2\int f^2 \d \mm h,$$
first for $f,\,g\in \dq(D)\cap C_0(E)$ and then 
for $f,\,g\in \dq(D)\cap C_b(E)$ (cf. [6; \S 3.2]).
Hence
$$D^{(c)}(fh, fh)
=\frac{1}{2}\int \d \mm{fh}
\le 2 D^{(c)}(f, f)+\int f^2 \d \mm h. \tag 3.3$$
On the other hand, since
$$|(fh)(x)- (fh)(y)|
\le |f(x)- f(y)| +I_{B\times A^c\cup A^c\times B}(x, y)
|f(1-h)(x)- f(1-h)(y)|,$$
we have
$$\align
&\int \jj [(fh)(x)- (fh)(y)]^2
\le 2\int \jj [f(x)- f(y)]^2\\
&\qd + 4\int_{B\times A^c}\jj
 [f(1-h)(x)- f(1-h)(y)]^2. \tag 3.4\endalign
$$
Thus, combining (3.1), (3.3), (3.4) with condition (2) together, we get
$$\align
D(fh, fh)&\le 2D(f, f)+\int f^2 \d \mm h
+2\int_{B\times A^c}\jj \big[f(1-h)(x)-f(1-h)(y)\big]^2\\
&\le 2 D(f, f) + 2c(h)\pi \big(f^2I_B\big)\\
&\le 2 D(f, f) + 2 \gz  c(h).	
\endalign$$
That is,
$$\align
D(f, f)&\ge \frac{1}{2}D(fh, fh) -\gz  c(h)
\ge \frac{1}{2}\lz_0\big(A^c\big)\pi \big(f^2h^2\big)-\gz  c(h)\\
&\ge \frac{1}{2}\lz_0\big(A^c\big)\pi \big(f^2I_{B^c}\big)-\gz  c(h)
= \frac{1}{2}\lz_0\big(A^c\big)(1-\gz ) -\gz   c(h).\tag 3.5
\endalign$$
Combining (3.2) with (3.5) together, we obtain
$$\align
D(f, f)&\ge \inf_{\gz \in [0, 1]}\max\bigg\{\frac{\lz_1(B)}{\pi(B)^2}
\big(\gz -\pi\big(B^c\big)\big),\,
\frac{1}{2}\lz_0\big(A^c\big)(1-\gz )-\gz  c(h)\bigg\}\\
&=\lz_1(B)\pi(B)^{-2}\big(\gz _0-\pi\big(B^c\big)\big). \tag 3.6
\endalign$$
The assertion of the theorem now follows from (3.6) and \lmm\;2.2.
\qed\enddemo
 
\thm\;3.1 is effective for diffusions was shown in [9] with a more direct
proof (in this case the Dirichlet form is explicit). We now apply the
theorem to jump processes.
 
\demo{Proof of \thm\;$1.4$} First, the topological assumptions of \thm\;3.1
are unnecessary in the present context. To see that condition (1) is
fulfilled, simply take $D_B$ to be the one defined by (1.13).
For condition (2), take $h=I_{A^c}$. Then
$$\int_{B\times A^c}\!\!\!\jj [(fI_A)(x)-(fI_A)(y)]^2
\!=\!\!\int_{A\times A^c}\!\!\!\jj f(x)^2
\!\le M_A\pi\big(f^2I_A\big)
\le M_A\pi\big(f^2I_B\big).
$$ 
This means that condition (2) holds with $c(h)=M_A$. We have thus proved
\thm\;1.4.\qed\enddemo
 
The application of \thm\;3.1 (or \thm\;1.4)
requires some estimates of $\lz_0(A^c)$ and
$\lz_1(B)$, which may be obtained from \thm s\;1.1---1.2. 
These estimates are usually in the qualitative sense good enough 
for $\lz_1(B)$, for which there are also quite a lot of publications 
including the authors' study in the past years. However, for $\lz_0(A^c)$,
the bound presented above may not be sharp enough, especially in the
unbounded situation. For this reason, we now introduce a different
result.
 
\proclaim{\thm\;3.2} Let $E$ be a metric space with Borel field $\eq$ and
let $(x_t)$ be a reversible right-continuous Markov process valued in $E$
with weak generator $\ooz$. Suppose that the corresponding Dirichlet form is
regular. Next, fix a closed set $B$. Suppose additionally that the
following conditions hold.
\roster
\item There exists a positive function $\fz$, $\fz|_B=0$ and
$$\sup_{B^c}\frac{1}{\fz} \ooz \fz =: \delta <0.$$
\item There exists a sequence of open sets $(E_n)$: 
$E_0\supset B$, $E_n\uparrow E$ such that $\fz$ is bounded below
on each $E_n\setminus B$ by a positive constant.
\item The first Dirichlet eigenfunction of $\ooz$ on each $E_n\setminus B$
is bounded above.
\endroster
Then we have $\lz_0\big(B^c\big)\ge  \delta $. In particular, for jump
processes, the condition ``$\fz|_B=0$'' given in (1) can be removed. 
\endproclaim  
 
Clearly, conditions (2) and (3) with compact $B$ are fulfilled for 
diffusions or Markov chains. Thus, the key condition here is the first
one.
 
\demo{\prf\; of \thm\;$3.2$} The last assertion follows by replacing
$\fz$ with $\fz I_{B^c}$. Indeed,
$$\align\ooz\big(\fz I_{B^c}\big)(x)&=
\int q(x, \d y)\big[\big(\fz I_{B^c}\big)(y)- 
\big(\fz I_{B^c}\big)(x)\big]\\
&\le \int q(x, \d y)\fz (y)- 
\big(\fz I_{B^c}\big)(x)
=\ooz\fz (x)\le -\delta\big(\fz I_{B^c}\big)(x) 
\qd\text{on}\qd B^c. \endalign$$
 
We are now going to prove the main assertion of the theorem. Set 
$\tz_B=\inf\{t\ge 0:\, x_t\in B\}$. Then, by condition (1) plus a truncating
argument if necessary, we get
$$\Bbb E^x\fz(x_{t\wedge \tz_B}) \le \fz (x),\qqd t\ge 0,\; x\notin B.$$
 
Next, let $u_n\,(\ge 0)$ be the first Dirichlet eigenfunction of $\ooz$
on $E_n\setminus B$. Set $\tz=\inf\{t\ge 0:\, x_t\notin E_n\setminus B\}$.
Then, by conditions (2) and (3), there exists $c_1>0$ such that
$u(x_{t\wedge\tz})\le c_1\fz (x_{t\wedge \tz_B})$ and so
$$u_n(x)e^{-\lz_0(E_n\setminus B)t}
=\Bbb E^x u_n(x_{t\wedge \tz})
\le c_1\Bbb E^x \fz(x_{t\wedge \tz_B})
\le c_1 \fz (x)e^{-\dz t},\qqd x\in E_n\setminus B.$$
This implies that $\lz_0(E_n\setminus B)\ge \dz$. Finally, because the
Dirichlet form is regular, it is easy to show that 
$\lz_0(B^c)=\lim_{n\to\infty} \lz_0(E_n\setminus B)$ and so the required
assertion follows. \qed
\enddemo
 
For the remainder of this section, we turn to study the upper bound of
$\lz_1$.
 
Let $(D, \dq(D))$ be a general conservative Dirichlet form and 
let $P(t, x, \d y)$ be the corresponding transition probability.
Fix $\fz\ge 0$. Suppose that $\fz\wedge n\in\dq(D)$ for every $n\ge 1$. 
Set $f_n=\exp[\vz (\fz \wedge n)/2]$.
Since the function $e^{\az x}$ is locally Lipschitz continuous and 
$\fz\wedge n$ is bounded, by the elementary spectral representation
theory, we have
$$\align
D(f_n, f_n)&=\lim_{t\to 0}\frac{1}{2t}\int \pi(\d x)P(t, x, \d y)
[f_n(x)-f_n(y)]^2\\
&\le \frac{\vz^2}{4} C(\fz, n)
\lim_{t\to 0}\frac{1}{2t}\int \pi(\d x)P(t, x, \d y)
[(\fz\wedge n)(x)-(\fz\wedge n)(y)]^2\\
&\le \frac{\vz^2}{4}C(\fz, n) D(\fz\wedge n, \fz\wedge n)<\infty,
\endalign$$
where $C(\fz, n)$ is the Lipschitz norm of $e^{\vz x/2}$ on the range of
$\fz\wedge n$. This leads us to introduce the following constant
$$\dz(\vz, \fz)=\vz^{-2}\sup_{n\ge 1} D(f_n, f_n)/\pi\big(f_n^2\big).$$
 
\proclaim{\thm\;3.3} Let $(D, \dq (D))$, $\fz$, $f_n$ and $\dz(\vz, \fz)$ 
be as above. Then, we have
$$\lz_1\le \sup\big\{\vz^2\dz(\vz, \fz):\, \pi\big(e^{\vz \fz}\big)
<\infty\big\}.$$
\endproclaim
 
\demo{\prf} We need to show that if $\pi\big(e^{\vz \fz}\big)
=\infty$, then $\lz_1\le \vz^2\dz(\vz, \fz)$. For $n\ge 1$, we have
$$\lz_1\le \frac{1}{2}
\frac{\int \jj [f_n(x)-f_n(y)]^2}{\pi\big(f_n^2\big)-\pi(f_n)^2}.
\tag 3.7$$
For every $m\ge 1$, choose $r_m>0$ such that $\pi(\fz\ge r_m)\le 1/m$.
Then
$$\pi \big(I_{[\fz\ge r_m]}f_n^2\big)^{1/2}
\ge \sqrt{m}\,\pi\big(I_{[\fz\ge r_m]}f_n\big)
\ge \sqrt{m}\,\pi(f_n) - \sqrt{m}\,e^{\vz r_m/2}.$$
Hence
$$\pi \big(f_n\big)^2
\le \Big[\sqrt{\pi\big(f_n^2\big)}\big/\sqrt{m} + e^{\vz r_m/2}
\Big]^2. \tag 3.8$$
 
On the other hand, by assumption, we have
$$D(f_n, f_n)\le \vz^2\dz(\vz, \fz)\pi\big(f_n^2\big).\tag 3.9$$
Noticing that $\pi\big(f_n^2\big)\uparrow \infty$, combining (3.9)
with (3.7) and (3.8) and then letting $n\uparrow \infty$, we obtain
$$\lz_1 \le \vz^2 \dz(\vz, \fz)/\big[1- m^{-1}\big].$$
The proof is completed by setting $m\uparrow \infty$.\qed
\enddemo
 
\demo{Proof of \thm\;$1.5$} It suffices to prove the first assertion 
because the remainder of the proof is similar. Let $f_n$ be given as
in \thm\;3.3. Note that by the Mean Value Theorem, 
$|e^A-e^B|\le |A-B| e^{A\vee B}=|A-B|(e^A\vee e^B)$ for all $A,\,B\ge 0$. 
Hence
$$\align
D(f_n, f_n)&=\frac{1}{2}\int \jj [f_n(x)-f_n(y)]^2\\
&\le \frac{\vz^2}{8}\int J^{(1)}(\d x, \d y)
[\fz(x)-\fz(y)]^2 r(x, y)\big[f_n(x)\vee f_n(y)\big]^2\\
&\le \frac{\vz^2}{4}\dz_2(\fz) \pi\big(f_n^2\big).
\endalign$$
The conclusion now follows from \thm\;3.3 with $\dz(\vz, \fz)=\frac{1}{4}
\dz_2(\fz)$. 
\qed\enddemo
 
\head 4. Existence of Spectral Gap for Markov Chains.\endhead
 
Usually, the power of a result for general jump processes should be 
justified by Markov chains.
 
Let $E$ be countable and $(q_{ij})$ be a regular and irreducible $Q$-matrix,
reversible with respect to $\pi=(\pi_i)$. As usual, let 
$q_i= \sum_{j\ne i} q_{ij}$. Assume that $q_i>0$ for all $i$ to rule out 
the reducible case. Then $K=0$ and 
$\ooz f(i)= \sum_{j\ne i} q_{ij}[f_j-f_i]$. 
The density of the symmetric measure with
respect to the counting measure becomes $J(i, j)=\pi_iq_{ij}\,(i\ne j)$. 
For simplicity, we consider
only two typical situations: $E=\Bbb Z_+$ or $E=\Bbb Z^d$
and take $r(i, j)=1/(q_i\vee q_j)$. Denote by
$|i|$ the $L^1$-norm, i.e., $|i|= \sum_{k=1}^d |i_k|$ for
$i=(i_1, \cdots, i_d)\in \Bbb Z^d$.
 
A combination of \thm\;1.2 and the next result provides us a 
simple condition for the existence of spectral gap 
for birth-death processes and the result seems to be new
from our knowledge even in such a simple situation (cf. [3]).
 
\proclaim{\thm\;4.1} Consider the birth-death process on $\Bbb Z_+$
with birth rates $(b_i)$ and death rates $(a_i)$.
\roster
\item Take $r_{ij}=(a_i+b_i)\lor (a_j+b_j)\,(i\ne j)$. 
Then $k^{(\az)'}>0\,($equivalently, $k^{(\az)}>0)$ iff 
there exists a constant $c>0$ such that
$$\frac{\pi_i a_i}{[(a_i+b_i)\vee (a_{i-1}+b_{i-1})]^{\az}}\ge
c \sum_{j\ge i}\pi_j, \qqd i\ge 1.\tag 4.1$$
Then, we indeed have $k^{(\az)'}\ge c$. Furthermore 
$$k^{(\az)}\ge \inf_{i\ge 1}\dfrac{\pz_i a_i}{
[(a_i+b_i)\vee (a_{i-1}+b_{i-1})]^{\az}(1-\pz_i) \sum_{j\ge i}\pz_j}.$$
\item Let $ \sum\limits_{i}\pz_i (a_i+b_i)\!<\!\infty$. 
Take $p_i=a_i+b_i$. Then we have 
$k_p'\!>\!0\, ($equivalently, $k_p\!>\!0)$ iff 
$\inf\limits_{i\ge 1}\dfrac{\pi_i a_i}{\sum_{j\ge i}\pi_j p_j}\!>\!0$ 
and moreover
$$k_p'\ge \inf_{i\ge 1}\dfrac{\pi_i a_i}{ \sum_{j\ge i}\pi_j p_j},\qqd
k_p\ge \inf_{i\ge 1}\dfrac{\pi_i a_i}{(1-\pz_i p_i/\bz_p)
\sum_{j\ge i}\pi_j p_j}.$$
\endroster 
\endproclaim
 
Roughly speaking, (4.1) holds if $\pi_j$ has exponential decay.
For polynomial decay, (4.1) can still be true when $\az=1/2$. See \xmp\;4.5.
 
\demo{\prf\; of \thm\;$4.1$} Here we prove part (1) only since the proof
of part (2) is similar.
 
a) Let $k^{(\az)}>0$. Take $A=I_i=\{i, i+1, \cdots\}$ for a fixed $i>0$
and
$$J^{(\az)}(i,j)=\frac{\pi_i q_{ij}}{[q_i\vee q_j]^{\az}}=
\cases
\dfrac{\pi_i a_i}{[(a_i+b_i)\vee (a_{i-1}+b_{i-1})]^{\az}}=:\pz_i\tilde a_i,
  \qd\text{if } j=i-1\\
\dfrac{\pi_ia_i}{[(a_i+b_i)\vee (a_{i+1}+b_{i+1})]^{\az}}=:\pz_i\tilde b_i,
  \qd\text{if }  j=i+1.
\endcases$$
Then
$$k^{(\az)'}\le k^{(\az)}\le \frac{J^{(\az)}(A\times A^c)}{\pi (A)\pi (A^c)}
=\frac{\pi_i\tilde a_i}{\big(\sum_{j\ge i}\pi_j\big) 
\big(\sum_{j< i}\pi_j\big)}
\le \frac{\pi_i\tilde a_i}{\pi_0 \sum_{j\ge i}\pi_j}.
$$
This proves the necessity of the condition.
 
b) Next, assume that the condition holds. Then for each $A$ with 
$\pi(A)\in (0,1)$, since the symmetry of $J^{(\az)}$, we may assume that
$0\notin A$. Set $i_0=\min A\ge 1$. Then, $A\subset I_{i_0}$,
$A^c\subset E\setminus \{i_0\}$ and so
$$\frac{J^{(\az)}(A\times A^c)}{\pi (A)\wedge \pi (A^c)}
\ge \frac{\pi_{i_0} \tilde a_{i_0}}{\sum_{j\ge i_0}\pi_j}
\ge c, \qqd
\frac{J^{(\az)}(A\times A^c)}{\pi (A)\pi (A^c)}
\ge \frac{\pi_{i_0} \tilde a_{i_0}}{(1-\pz_{i_0})\sum_{j\ge i_0}\pi_j}.$$
Because $A$ is arbitrary, we obtain the required assertions.
\qed\enddemo
 
\proclaim{\thm\;4.2} Let $E=\Bbb Z_+$. Suppose that 
$(q_{ij})$ has finite range $R$, i.e., 
$q_{ij}=0$ whenever $|i-j|>R$. Then, we have
$\lz_1>0$ provided
$$\varlimsup_{i\to\infty}\sum_{j}\frac{q_{ij}}{\sqrt{q_i\vee q_j}} (j-i)<0.$$
\endproclaim
 
\demo{\prf} Simply take $\fz_i=i+1$ and $B=\{0, 1, \cdots, n\}$ for large
$n$ in \thm\;1.3 and then apply \thm\;1.2.\qed\enddemo
 
Similarly, we have the following result.
 
\proclaim{\thm\;4.3} Let $E=\Bbb Z^d$. Suppose that 
$(q_{ij})$ has finite range $R$. Then, we have
$\lz_1>0$ provided
$$\varlimsup_{|i|\to\infty}
\sum_{j}\frac{q_{ij}}{\sqrt{q_i\vee q_j}}\sum_{k=1}^d \big[|j_k|\vee  R-
|i_k|\vee  R\big]<0.$$
\endproclaim
 
\demo{\prf} Take $\fz_i= \sum_{k=1}^d |i_k|\vee  R+1$ in
\thm\;1.3 and then apply \thm\;1.2.\qed\enddemo
 
\proclaim{\thm\;4.4} Let $E=\Bbb Z^d$. 
If there exists a positive function $\fz$ such that
$$\varlimsup_{|i|\to\infty}\frac{1}{\fz}\ooz \fz<0,$$
then $\lz_1>0$.
\endproclaim
 
\demo{\prf} Apply \thm\;1.2, \thm\;3.2 and then \thm\;1.4 to the finite sets
$\{i: |i|\le n\}$.
\qed\enddemo
 
The following example, taken from [3], is especially rare and 
interesting since it 
exhibits the critical phenomena for the existence of spectral gap. 
It is now used to justify the power of our results and we should see soon 
what will happen.
Similar example for diffusion was given in [4] and [9].
 
\proclaim{\xmp\;4.5} Let $E=\Bbb Z_+$ and $a_i=b_i=i^\gz\,(i\ge 1)$ for
some $\gz>0$, $a_0=0 $ and $b_0=1$. Then $\lz_1>0$ iff $\gz\ge 2$.
\endproclaim 
 
\demo{\prf} a) By part (1) of \thm\;4.1, we have $k^{(1/2)}>0$ iff 
$\gz \ge 2$. Thus, by \thm\;1.2, we have $\lz_1>0$ for all $\gz\ge 2$. 
 
b) Applying \thm\;1.5 to $\fz_i=1+i^{1-\gz/2}$, it follows that $\lz_1=0$ 
for all $\gz\in (1, 2)$. 
 
c) The conditions of \thm\;4.2 hold whenever $\gz\ge 2$. Hence $\lz_1>0$
for all $\gz\ge 2$.
 
d) Next, taking $\fz_i=\sqrt{i}\;(i\ge 1)$, we see that 
$\ooz \fz(i)/\fz_i=-\frac{1}{2\gz (\gz-1)}i^{\gz -2}+O\big(i^{\gz -3}\big)$.
Then
$$\lim_{i\to\infty}\frac{1}{\fz_i}\ooz \fz(i)=
\cases
-\infty, \qd & \gz>2\\
-\frac{1}{4}, \qd & \gz=2.
\endcases$$ 
By \thm\;4.4, we have $\lz_1>0$ for all $\gz\ge 2$.
 
On the other hand, take $f_n(i)=i^{\frac{\gz-1}{2}}\wedge 
n^{\frac{\gz-1}{2}}$ and $A=\{0\}$. Then
$$\align
\lz_0(A^c)&\le
\varliminf_{n\to\infty}\frac{ \sum_{i,j\ge 0} \pi_iq_{ij}[f_n(j)-f_n(i)]^2}
{2 \sum_{i\ge 0} \pi_i f_n(i)^2}
=\varliminf_{n\to\infty}\frac{ \sum_{i\ge 0} \pi_i
q_{i, i+1}[f_n(i+1)-f_n(i)]^2}
{2 \sum_{i\ge 0} \pi_i f_n(i)^2}\\
&\le\varliminf_{n\to\infty}\frac{1+ (\gz-1)^2  \sum_{i=1}^n i^{\gz-3}}
{\sum_{i=1}^n i^{-1}}=0, \qqd 1< \gz <2.
\endalign$$
By \thm\;1.4, we get $\lz_1\le \lz_0(A^c)/\pi (A)=0$. The case that $\gz
\le 1$ can be ignored since then the chain is not positive recurrent.
\qed\enddemo
 
Thus, we have seen that all the results presented in the paper, except
\thm\;2.5 which does not work for this example, are qualitatively sharp 
for this example since every one covers the required region 
and there is no gap left.
Finally, taking $\az=0$ in part (1) of \thm\;4.1, we obtain
$k\ge \big( \sum_{i=1}^\infty i^{-\gz}\big)^{-1}>0$ for all $\gz>1$.
In other words, we have $k>0$ but $\lz_1=0$ for all $\gz\in (1, 2)$.
Therefore, the condition ``$k>0$'' is not but ``$k^{(1/2)}>0$'' is
sufficient for $\lz_1>0$. 
 
The next two examples show that the two approaches used in the paper
for the Cheeger's inequalities may all attain sharp estimates but they 
are not comparable (remember that \thm\;2.5 is not suitable for
\xmp\;4.5). We mention that as far as we know, no optimal estimate
provided by the Cheeger's technique ever appeared before.
 
\proclaim{\xmp\;4.6} Let $E=\Bbb Z_+$ and take $a_i\equiv a>0$ and
$b_i\equiv b>0$. Then, both \thm\;1.2 and \thm\;2.5 are sharp.
\endproclaim
 
\demo{\prf} This is a standard example which is often used to justify
the power of a method. It is well known that 
$\lz_1= \big(\sqrt{a}-\sqrt{b}\,\big)^2$ (cf. [2; \xmp\;9.22] and [3]).
 
a) By part (1) of \thm\;4.1, we have  
$${k^{(\az)}}'\ge \inf_{i\ge 1}
\frac{\pi_i a_i}{(a+b)^{\az} \sum_{j\ge i} \pi_j}
=\frac{a-b}{(a+b)^{\az}}.$$
Then, by \thm\;1.2, we get $\lz_1\ge \big(\sqrt{a}-\sqrt{b}\,\big)^2$.
 
b) Take $p_i\equiv a+b$. Then by part (2) of \thm\;4.1,
$$k'_p\ge \inf_{i\ge 1}\frac{\pi_i a_i}{\sum_{j\ge i} \pi_j p_j}
=\frac{a-b}{a+b}. $$
The same estimate as in a) now follows from \thm\;2.5.\qed\enddemo
 
\proclaim{\xmp\;4.7} Let $E=\Bbb Z_+$ and take $q_{0k}=\bz_k>0$
(be careful to distinguish the sequence $(\bz_k)$ and the constant $\bz_p$), 
$q_{k0}=1/2$ for $k\ge 1$ and $q_{ij}=0$ for all other $i\ne j$. Assume that 
$q_0= \sum_{k\ge 1}\bz_k<\infty$.
Then, \thm\;2.5 is sharp for all $q_0$ but \thm\;1.2 is sharp only for
$q_0\le 1/2$.
\endproclaim
 
\demo{\prf} From $\pz_0q_{0k}=\pi_k q_{k0}$, it follows that 
$\pz_k=2\pz_0\bz_k$, $k\ge 1$ and $\pz_0=(1+2q_0)^{-1}$. An interesting
point of the example is that the decay of $ \sum_{j\ge i} \pi_j$ as
$i\to\infty$ can be arbitrary slow, not necessarily exponential. The last
condition is necessary for $\lz_1>0$ for the birth-death processes
with rates bounded below (by a positive constant) and above 
(cf. [2; \crl\;9.19 (4)]).
 
a) Take $p_i=q_i\vee (1/2)$, then $\az_p=1/2$.
Without loss of generality, assume that $0\notin A$. Then
$$\align
\frac{1}{\bz_p}\cdot\frac{J (A\times A^c)}{\pi_p(A)\wedge \pi_p(A^c)}
&=\frac{ \sum_{i\in A} \pz_i q_{i0}}
{\Big( \sum_{i\in A}2\pi_0\bz_ip_i\Big)
\wedge\Big(\pz_0 p_0+ \sum_{i\notin A,\,i\ne 0 }2\pz_0\bz_ip_i\Big)}\\
&=\frac{ \sum_{i\in A} \bz_i }
{\Big( \sum_{i\in A}2 \bz_i p_i\Big)
\wedge\Big(p_0+ \sum_{i\notin A,\,i\ne 0 }2\bz_i p_i\Big)}\\
&=\frac{ \sum_{i\in A} \bz_i }
{\Big( \sum_{i\in A} \bz_i \Big)
\wedge\Big(p_0+ \sum_{i\notin A,\,i\ne 0 }\bz_i \Big)}\ge 1.
\endalign$$
This gives us $k_p'\ge 1$ and hence by \thm\;2.5,
$$\lz_1\ge \az_p\Big(1-\sqrt{1-{k_p'}^2}\,\Big)\ge 1/2.$$
Actually, every equality in the last line must hold.
 
b) Again, assume that $0\notin A$. Then
$$\align
\frac{J^{(\az)}(A\times A^c)}{\pi(A)\wedge \pi(A^c)}
&=\frac{ \sum_{i\in A} \pz_i q_{i0}(q_i\vee q_0)^{-\az}}
{\Big( \sum_{i\in A}2\pz_0\bz_i\Big)
\wedge\Big(\pi_0+ \sum_{i\notin A,\,i\ne 0 }2\pz_0\bz_i\Big)}\\
&=\frac{1}{2}\cdot\frac{ \sum_{i\in A} 2\bz_i }
{\big(\frac{1}{2}\vee q_0\big)^{\az}
\Big[\Big( \sum_{i\in A}2 \bz_i\Big)
\wedge\Big(1+ \sum_{i\notin A,\,i\ne 0 }2\bz_i\Big)\Big]}\\
&=\frac{1}{2}\cdot \frac{1}{\big(\frac{1}{2}\vee q_0\big)^{\az}}\cdot
\frac{ \sum_{i\in A} \bz_i }
{\Big( \sum_{i\in A} \bz_i\Big)
\wedge\Big(1/2+ \sum_{i\notin A,\,i\ne 0 }\bz_i\Big)}\\
&=\frac{1}{2}\cdot \frac{1}{\big(\frac{1}{2}\vee q_0\big)^{\az}}\cdot
\frac{1}
{1\wedge \Big[\Big(1/2+ \sum_{i\notin A,\,i\ne 0 }\bz_i\Big)
\Big/\sum_{i\in A} \bz_i\Big]}.
\endalign$$
Because $\Big(1/2+ \sum_{i\notin A,\,i\ne 0 }\bz_i\Big)
\Big/\sum_{i\in A} \bz_i$ decreases when $A$ increases, by setting
$A=\{i\}$ for a large enough $i\ne 0$, it follows that
$${k^{(\az)}}'=\inf_{A:\, 0\notin A}
\frac{J^{(\az)}(A\times A^c)}{\pi(A)\wedge \pi(A^c)}
=\frac{1}{2}\bigg(\frac{1}{2}\vee q_0\bigg)^{-\az}.$$
By \thm\;1.2, we get
$$\lz_1\ge \frac{1}{2}\Big\{ 1\vee (2q_0)+
\sqrt{\big(1\vee (2q_0)\big)^2-1}\,\Big\}^{-1}.$$
Thus, the lower bound is equal to $1/2=\lz_1$ iff $q_0\le 1/2$.
\qed\enddemo
 
The following counterexample shows the limitation of the Cheeger's 
inequalities. Of course, the example can be easily handled with the help of
some comparison technique. However, this suggests us that sometimes it
is necessary to examine a model carefully before applying the
inequalities.  
 
\proclaim{\xmp\;4.8} Consider the birth-death process with
$a_{2i-1}=(2i-1)^2$, $a_{2i}=(2i)^4$ and
$b_i=a_i$ for all $i\ge 1$. Then, we have $k^{(1/2)'}=0$ and so
\thm\;1.2 is not applicable.
\endproclaim
 
\demo{\prf} First, applying \thm\;4.4 to $\fz_i=\sqrt{i}$ or comparing the 
chain with the one with rates $a_i=b_i=(2i)^2$, one sees that $\lz_1>0$.
Next, because $\mz_i=1/a_i$ (and hence $\pz_i=\mz_i/Z$,
where $Z$ is the normalizing constant), we have
$ \sum_{j\ge i}\mz_j=O(i^{-1})$. However, $\sqrt{a_i\vee a_{i-1}}=O(i^2)$.
Hence $\sup_{i\ge 1} \sqrt{a_i\vee a_{i-1}}\,\sum_{j\ge i}\mz_j$
\llnbr$=\infty$.
This gives us $k^{(1/2)'}=0$ by part (1) of \thm\;4.1.
 
Note that the choice $r_{ij}=q_i\vee q_j\,(i\ne j)$ is usually not optimal
in the sense for which (1.8) often becomes inequality rather than equality.
However, the improvement provided by an optimal $r_{ij}$ is still not
enough to cover this example and so the problem is really due to the
limitation of the technique.
\qed\enddemo

\bigskip
 
\head{5. Appendix. {\prf\; of \lmm\;$2.1$}, for referee's reference}
\endhead
 
Because $\az \ge 0$ is fixed, we can
omit the superscript ``$(\az)$'' everywhere in the proof.
Denote by $\tilde h$, $\tilde k$ and $\tilde k'$ the 
right-hand sides of the above quantities. 
By taking $f=I_A$, we obtain $h\ge \tilde h$, $k\ge \tilde k$
and $k'\ge \tilde k'$. We now prove the reverse inequalities.
 
a) For any $f\ge 0$ with $\pi (f)=1$, let $A_{\gz}=\{f> \gz\}$, $\gz\ge 0$.
By the symmetry of $J$, we have
$$\align
&\frac{1}{2}\int \jj |f(y)-f(x)|+K(f)
= \int_{\{f(x)>f(y)\}} \jj [f(x)-f(y)]+K(f)\\
&=\int_0^\infty\d \gz\Big\{J\big(\{f(x)>\gz\ge f(y)\}\big)+
K\big(\{f>\gz\}\big)\Big\}\\
&=\int_0^\infty\big [J\big(A_\gz\times A_\gz^c\big)+
K(A_\gz)\big]\d \gz
\ge h\int_0^\infty \pi\big(A_\gz\big)\d \gz=h \pi (f)=h.\endalign
$$
Hence $\tilde h\ge h$.
 
b) For any $f\in L_+^1(\pi)$ with $\int\pi(\d x)\pi(\d y)
|f(x)-f(y)|=1$, by a), we have
$$\align
&\int\jj |f(x)-f(y)|=2\int\d \gz J(A_{\gz}\times A_{\gz}^{c})\\
&\ge k \int_0^\infty\d \gz \;(\pi\times \pi)(A_{\gz}\times A_{\gz}^c)
=k\int_0^\infty \pi(\d x)\pi (\d y)|f(x)-f(y)|=k.
\endalign$$
This proves the first equality of $k^{(\az)}$.
 
Next, we show that
$$\int |f-\pi (f)|\d \pi=\sup_{g:\,\pi (g)=0,\,\inf_{c\in \Bbb R}
\|g-c\|_\infty\le 1}\int fg\d \pi,\tag 5.1$$
where $\|\cdot\|_p$ denotes the $L^p$-norm.
First, let $\pi(g)=0$ with $\inf_{c\in \Bbb R}\|g-c\|_\infty\le 1$. Then,
because $\pi(g)=0$ and $\pi (f-\pi (f))=0$, we have
$$\int fg\d \pi=\int(f-\pi (f))g\d\pi
=\int (f-\pi(f))(g-c)\d\pi$$
for all $c\in \Bbb R$. Hence, by H\"older inequality, we have
$$\bigg|\int fg\d\pi\bigg|\le \|f-\pi(f)\|_1\|g-c\|_{\infty}$$
for all $c$. This gives us 
$$\bigg|\int fg\d\pi\bigg|\le \|f-\pi(f)\|_1\,\inf_c\|g-c\|_{\infty}
\le \|f-\pi(f)\|_1.$$
On the other hand, for a given $f\in L^1(\pi)$, set $A_f^+=\{f\ge\pi(f)\}$
and $A_f^-=\{f<\pi(f)\}$. Take 
$g_0=I_{A_f^+}-I_{A_f^-}-\pi(A_f^+)+\pi(A_f^-)$.
Then, $g_0\in L^\infty(\pi)$ and $\pi (g_0)=0$. Finally, 
take $c_0=1-2\pi(A_f^+)$.
Then, it is easy to check that $\inf_c\|g_0-c\|_\infty=\|g_0-c_0\|_\infty=1$.
Therefore, we have $\int fg_0\d \pi=\int|f-\pi(f)|\d\pi$ as required.
 
We now prove the second equality of $k^{(\az)}$. 
Let $f\ge 0$ and set $A_{\gz}=\{f\ge \gz\}$. Again, by using a) and (5.1), 
we have
$$\align
&\int\jj |f(y)-f(x)|\ge 
2 k\int_0^\infty \d\gz \pi(A_{\gz}) \pi(A_{\gz}^c)
= k\int_0^\infty \d\gz \int\big|I_{A_{\gz}}-\pi(A_{\gz})\big|\d\pi\\
&\qd= k\int_0^\infty\d \gz\sup_{g:\,\pi (g)=0,\,
\inf_{c\in \Bbb R}\|g-c\|_\infty\le 1}
\int I_{A_{\gz}}g\d \pi\\
&\qd\ge k \sup_{g:\,\pi (g)=0,\, \inf_{c\in \Bbb R}\|g-c\|_\infty\le 1} 
\int_0^\infty\d \gz \int I_{A_{\gz}}g\d \pi\\
&\qd= k \sup_{g:\,\pi (g)=0,\, \inf_{c\in \Bbb R}\|g-c\|_\infty\le 1} 
\int f g\d \pi= k\int |f-\pz(f)|\d \pz.
\endalign$$
Therefore, we obtain $\tilde k\ge k$.
 
c) Choose $c_0\in \Bbb R$ such that $\pi (f<c_0)$, $\pi (f>c_0)\le 1/2$.
Let $f_{\pm}=(f-c_0)^{\pm}$. Then we have $f_++f_-=|f-c_0|$ and
$\pi (|f-c_0|)=\min_c\pi (|f-c|)$. For any $\gz\ge 0$, define
$A_\gz^\pm=\{f_\pm>\gz\}$. We have
$$\align
&\frac{1}{2}\int \jj |f(y)-f(x)|=\frac{1}{2}
 \int \jj \big[|f_+(y)-f_+(x)|+ |f_-(y)-f_-(x)|\big]\\
&=\int_0^\infty\big [J\big(A_\gz^+\times A_\gz^{+c}\big)+
J\big(A_\gz^-\times A_\gz^{-c}\big)\big]\d \gz\\
&\ge k'\int_0^\infty \big[\pi\big(A_\gz^+\big)+ \pi\big(A_\gz^-\big)\big]
\d \gz=k' \pi (f_++f_-)=k'\pi (|f-c_0|)=
k'\min_c \pi (|f-c|).\endalign
$$
This implies that $\tilde k'\ge k'$.\qed

\demo{Acknowledgement} The second named author would like to thank MSRI for
a happy stay.\enddemo
 

\Refs
 
\widestnumber\no {100}
 
\ref\key 1
\by Cheeger, J. (1970)
\paper A lower bound for the smallest eigenvalue of the Laplacian
\jour Problems in analysis, a symposium in honor of S. Bochner
195--199, Princeton U. Press, Princeton\endref
 
\ref\key 2 \by Chen, M. F. (1992)
\book From Markov Chains to Non-Equilibrium Particle Systems
\publ  Singapore, World Scientific \endref
 
\ref\key 3\by Chen, M. F. (1996)
\paper Estimation of spectral gap for Markov chains
\jour Acta Math. Sin. New Ser. 12:4, 337--360\endref
 
\ref\key 4
\by Chen, M. F. and Wang, F. Y. (1997)
\paper Estimation of spectral gap for elliptic operators
\jour Trans. Amer. Math. Soc. 349, 1209--1237\endref
 
\ref\key 5
\by Chung, F. R. K. (1997)
\book Spectral Graph Theory
\publ CBMS, {\bf 92}, AMS, Providence, Rhode Island
\endref
 
\ref\key 6
\by Fukushima, M., Oshima, Y. and Takeda, M. (1994)
\book Dirichlet Forms and Symmetric Markov Processes
\publ Walter de Gruyter \& Co
\endref
 
\ref\key 7
\by Lawler, G. F. and Sokal, A. D. (1988)
\paper Bounds on the $L^2$ spectrum for Markov chain and Markov 
processes: a generalization of Cheeger's inequality
\jour Trans. Amer. Math. Soc. {\bf 309}, 557--580\endref
 
\ref\key 8
\by Saloff-Coste, L. (1997)
\paper Lectures on Finite Markov Chains
\jour LNM {\bf 1665}, Springer-Verlag\endref
 
\ref\key 9
\by Wang, F. Y. (1997)
\paper Existence of spectral gap for diffusion processes
\jour preprint\endref
\endRefs
\enddocument